\documentclass[a4paper,10pt,fleqn]{article}

\usepackage{theorem}
\usepackage{amssymb}
\usepackage{latexsym}

\theoremstyle{break}
\newtheorem{thm}{Theorem}[section]
\newtheorem{prop}{Proposition}[section]
\newtheorem{lem}{Lemma}[section]
\newtheorem{de}{Definition}[section]

\title{\ 
 M\"obius numbers of some modified generalized noncrossing partitions}
\author{Masaya Tomie}
\date{tomie@math.tsukuba.ac.jp}

\begin{document}

\maketitle

In this paper we will give a M\"obius number of 
$NC^{k}(W) \setminus \bf{mins} \cup \{  \widehat{0} \}$ for a Coxeter group 
$W$ which contains an affirmative answer for the conjecture 3.7.9 in 
Armstrong's paper 
[ Generalized noncrossing partitions and combinatorics of Coxeter groups. 
arXiv:math/0611106].

\section{Introduction}

In this paper we will prove the following theorem which is conjectured in 
 \cite{armstrong}.

\begin{thm}\label{theorem}

For each finite Coxeter group $(W,S)$ with $|S| = n$ and 
for all positive integers  $k$, the M\"obius number of  
$NC^{k}(W) \setminus \bf{mins} \cup \{  \widehat{0} \}$ equals to  
$(-1)^{n} \Bigl(  Cat_{+}^{(k)}(W) - Cat_{+}^{(k-1)}(W)  \Bigr)$.

\end{thm}

Our method is using 
 the EL-labeling of $NC_{(k)}(W)$ introduced by Armstrong and Thomas 
 \cite{armstrong}. 
If we  give an EL-labeling for $NC(W)$  for any complex reflction group $W$, 
 then we can state our Theorem {\rmfamily \ref{theorem}} in the case of 
  any well-generated complex reflection group. 
 But it may be difficult to give a uniform proof 
   because Athanasiadis, Brady and Watt gave an EL-labeling for 
   $NC(W)$ using some properties of the root system derived from  
   a real reflection group $W$ \cite{abw}.
In \cite{armstrong-krattenthaler} they proved this result by counting 
the multichains of $NC_{(k)}(W)$. 
Moreover they proved in the case of well-generated complex reflection 
groups.
 Our approach is independent to theirs.  It is surprising for us that   
their paper \cite{armstrong-krattenthaler} 
appeared in arXiv when we were typing this paper.

\section{Preliminaries}

\subsection{generalized noncrossing partition}
 
In this paper we put $(W,S)$ a Coxeter group $W$ with a set of generators 
$S$ where $S=n$. Basic properties of  Coxeter groups is 
introduced in \cite{humphleys}.
We put $T:= \{ wsw^{-1} \ | \ s \in S, w \in W \}$ the 
cojugate closure  of the set of generators $S$.
Let $l_{T} \ : \ W \longrightarrow \mathbb{Z}$ denote the word length on $W$ 
with respect to the set $T$. 
We call $l_{T}$ the absolute length on $W$. 
Then the absolute length  naturally 
induces a partial order on $W$ as following: 
$\pi \le \sigma$ if $l_{T}(\sigma) = l_{T}(\pi) + l_{T}(\pi^{-1} \sigma)$.
We call it the absolute order on $W$.
We fix a Coxeter element $\gamma \in W$ and call the poset $[e,\gamma]$ 
with the absolute order $NC(W)$. 
Next we put $NC^{(k)}(W) := \{ (\pi_{1}, \ldots ,\pi_{k}) \ | \ \pi_{i} \in 
NC(W) \ {\rm for \ } 1 \le i \le k \  {\rm with } \ \pi_{1} 
\le \pi_{2} \le \cdots \le \pi_{k} \le \gamma \}$ and 
$NC_{(k)}(W) := \{ (\delta_{1}, \ldots ,\delta_{k}) \ | \ \delta_{i} \in 
NC(W) \ {\rm for } \ 1 \le i \le k \ {\rm with } \ 
l(\delta_{1} \cdots \delta_{i} ) = l(\delta_{1}) + \cdots l(\delta_{i}) 
\ {\rm for } \ 1 \le i \le k \}$. In \cite{armstrong} Armstrong introduced the 
 order structure of them is as follows:
 For $(\pi)_{k}^{(1)} := (\pi_{1}^{(1)}, \ldots ,\pi_{k}^{(1)}) \ {\rm and} \ 
 (\pi)_{k}^{(2)} := (\pi_{1}^{(2)}, \ldots ,\pi_{k}^{(2)}) \in 
 NC^{(k)}(W)$, $(\pi)_{k}^{(1)} \le (\pi)_{k}^{(2)}$ if 
 $ (\pi_{i}^{(2)})^{-1} (\pi_{i+1}^{(2)}) 
 \le 
 (\pi_{i}^{(1)})^{-1} (\pi_{i+1}^{(1)})  $
  for $1 \le i \le k$ where $\pi_{k+1}^{(1)} = \pi_{k+1}^{(2)} = \gamma$.
  
  For $(\delta)_{k}^{(1)} := (\delta_{1}^{(1)}, \ldots ,\pi_{k}^{(1)}) \ 
  {\rm and} \ 
 (\delta)_{k}^{(2)} := (\delta_{1}^{(2)}, \ldots ,\delta_{k}^{(2)}) \in 
 NC_{(k)}(W)$, 
 $(\delta)_{k}^{(1)} \le (\delta)_{k}^{(2)}  \ {\rm if } \ 
 \delta_{i}^{(1)} \le \delta_{i}^{(2)}$ for $1 \le i \le k$.
 It is easy to see that the poset $NC_{(k)}(W)$ is the dual poset of 
 $NC^{(k)}(W)$ (for more information, see \cite{armstrong}).

\subsection{EL-shellability}

Let $(P, \preceq )$ be a finite poset. Assume that $P$ is bounded, meaning 
that $P$ has a minimum element and a maximum element, denoted 
$\widehat{0}$ and $\widehat{1}$ respectively, and that it is graded, 
meaning that all maximal chains in $P$ have the same length. This length is 
called the $rank$ of $P$ and denoted rank$(P)$.
Let $\epsilon (P)$ be the set of covering relations of $P$, meaning pairs 
$(x,y)$ of elements of $P$ such that $x \prec y$ in $P$.
Let $\Lambda$ be a totally ordered set. An $edge \ labeling $ of $P$ with 
label set $\Lambda$ is a map $\lambda \ : \ \epsilon(P) \longrightarrow 
\Lambda$. 
Let $c$ be an unrefinable chain $x_{0} \prec x_{1} \prec \cdots \prec x_{r}$ 
of elements of $P$ so that $(x_{i-1},x_{i}) \in \epsilon(P)  $
 for all $1 \le i \le r$. 
 We let $\lambda(c) = (\lambda(x_{0},x_{1}), \lambda(x_{1},x_{2}), \cdots 
 \lambda(x_{r-1},x_{r}))$ be the label of $c$ with respect to $\lambda$ and 
 call $c$ $rising$ and $falling$ with respect to  $\lambda$ if the entries 
 of $\lambda(c)$ strictly increase or weakly decrease, respectively, in the 
 total order of $\Lambda$. 
 We say that $c$ is $lexicographically \ smaller$ than an unrefinable chain 
 $\acute{c}$ in $P$ with respect to $\lambda$ if $\lambda(c)$ proceeds 
 $\lambda(\acute{c})$ in the lexicographic order induced by the total order 
 of $\Lambda$ \cite{abw}.

 \begin{de}[\cite{bjo}]
 
 An edge labeling $\lambda$ of $P$ is called an EL-labeling if for every 
 nonsingleton interval $[u,v]$ in $P$
 
 (1) there is a unique rising maximal chain in $[u,v]$ and 
 
 (2) this chain is lexicographically smallest among all maximal chains in 
 $[u,v]$
 
 with respect to $\lambda$.

 \end{de}
 
 The poset $P$ is called EL-shellable if it has EL-labeling for some 
 label set $\Lambda$.
 
 For a graded and bounded poset $(P,\preceq)$, we denote by $\mu(P)$ 
 the M\"obius number of $P$.
 If $P$ is EL-shellable the M\"obius number of 
 $P$ is the number of falling maximal chains of $P$ up to 
 sign $(-1)^{rank(P)}$ \cite{stan}.

\section{Main result}

In this section we will prove the following Theorem.

\begin{thm}

For each finite Coxeter group $(W,S)$ with $|S| = n$ and 
for all positive integers  $k$, we have 
$\mu(NC^{k}(W) \setminus \bf{mins} \cup \{  \widehat{0} \})= $ 
$(-1)^{n} \Bigl(  Cat_{+}^{(k)}(W) - Cat_{+}^{(k-1)}(W)  \Bigr)$.

\end{thm}

It is easy to see the following Lemma.

\begin{lem}

Let $P$ be a graded  poset with a minimum element $\widehat{0}$. 
We put $\bf{maxs(P)}$ the set of maximal elements of $P$. 
Then the poset $P \setminus \bf{maxs(P)}$ is also graded.
We denote by $\mu(P \setminus \bf{maxs(P)} \cup \{\widehat{1}\})$ 
 the M\"obius number of $P \setminus \bf{maxs(P)} \cup \{\widehat{1}\}$.
 Then we have 
 $\mu(P \setminus \bf{maxs(P)} \cup \{\widehat{1}\}) = 
 \mu(P \cup \{\widehat{1}\}) 
 + \Sigma_{x \in \bf{maxs(P)}} \mu([\widehat{0}
, x])$

\end{lem}

For $k \in \mathbb{N}$ and an arbitrary finite Coxeter group $(W,S)$ 
we consider the poset $NC_{(k)}(W)$ which is the dual poset of 
$NC^{(k)}(W)$. 
We put ${\bf maxs}$ as a set of maximal elements of $NC_{(k)}(W)$.
To show our Theorem, it is sufficient to prove 
$\mu (NC_{(k)}(W) \setminus {\bf maxs} \cup \{ \widehat{1}    \}) = 
(-1)^{n} \Bigl(  Cat_{+}^{(k)}(W) - Cat_{+}^{(k-1)}(W)  \Bigr) $

\vspace{5mm}

In \cite{armstrong} Armstrong and Thomas gave an EL-shelling of 
$ NC_{k}(W) \cup \{ \widehat{1}   \}$. 
We will explain their method briefly.
We put $T$ the set of reflections of $W$. 
Recall that the edges in the Hasse diagram of $NC(W)$ are naturally labelled 
by reflections $T$. Athanasiadis, Brady and Watt defined a total order on 
the set $T$ such that the natural edge-labelling by $T$ becomes an 
EL-shelling of the poset $NC(W)$. 
We put the EL-labeling $\lambda \ : \ \epsilon(NC(W)) \longrightarrow T$.
In \cite{abw} they called 
the total order on $T$ the {\bf ABW} order. 
They put $T:= \{ t_{1} , \cdots t_{N} \}$ with 
the {\bf ABW} order $t_{1} < t_{2} < \cdots < t_{N}$.
Recall that $NC(W^{k})$ is edge-lebelled by the set of reflections 
$T^{k} := \{ t_{i.j} = (1, 1, \cdots ,t_{i,j}, \cdots , 1 ) :  1 \le i,j \le 
N \}$ where $t_{j}$ occurs in the $i$-th entry of $t_{i,j}$.
Then they defined the {\bf lex ABW order} on $T^{k}$ as 
$t_{1,1} < t_{1,2} < \cdots < t_{1,N} < t_{2,1} < t_{2,2} < \cdots 
t_{2,N} <  \cdots t_{k,1} < t_{k,2} < \cdots t_{k,N}$.
This induces an EL-shelling of $NC(W^{k})$.
Now recall that $NC_{k}(W)$ is an order ideal in $NC(W^{k})$, so the 
{\bf lex ABW order} on $T^{k}$ restricts to an EL-labelling of the Hasse 
diagram of $NC_{k}(W)$. 
They considered the set $T^{k} \cup \{ \theta \}$ with 
$t_{1,1} < t_{1,2} < \cdots < t_{1,N} < \lambda < t_{2,1} < t_{2,2} < \cdots 
t_{2,N} <  \cdots t_{k,1} < t_{k,2} < \cdots t_{k,N}$.
For $x \in {\bf maxs } $ they put $\lambda(x, \widehat{1}) := \lambda$, where 
$\lambda(x, \widehat{1})$ is the edge from $x$ to $\widehat{1}$. They showed 
that the labeling as above induces  an EL-shelling of 
$NC_{k}(W) \cup \{\widehat{1}\}$.
Now we put their EL-labeling 
$\widehat{\lambda} \ : \ \epsilon(NC_{k}(W) \cup \{\widehat{1}\}) 
\longrightarrow T^{k} \cup \{ \theta \}$.

We have 

$\mu(NC^{(k)}(W) \setminus {\bf mins} \cup \{\widehat{0}\})$ 
$= \mu(NC_{(k)}(W) \setminus {\bf maxs} \cup \{\widehat{1}\})$

$\Sigma_{x \in {\rm maxis}} \mu(\widehat{0},x) + 
\mu(NC_{(k)}(W) \cup \{\widehat{1}\})$

$= \Sigma_{x \in {\rm maxs}} 
\mu(\widehat{0},x) + (-1)^{n-1} $Cat$_{+}^{(k-1)}(W)$.

It is sufficient to show 
$\Sigma_{x \in {\rm maxs}} \mu(\widehat{0},x) = (-1)^{n} $Cat$_{+}^{(k)}(W)$ to 
prove Theorem {\rmfamily \ref{theorem}}.
First we consider the EL-shelling introduced by Armstrong and Thomas.
Recall that $\mu(NC_{k}(W) \cup \{\widehat{1}\}) = (-1)^{n-1} \times$ 
 the number of the falling  maximal chains of 
$NC_{k}(W) \cup \{\widehat{1}\}$ with respect to $\widehat{\lambda}$.

Now let $c$ be an unrefinable  chain 
$( e, \cdots ,e) \prec \cdots \prec (\delta_{1}, \cdots \delta_{k})
\prec \widehat{1} $ of elements of $NC_{k}(W) \cup \{\widehat{1}\}$.
If $c$ is a falling maximal chain with respect to $\widehat{\lambda}$, 
we must have $\delta_{1} = e$  because 
$\widehat{\lambda}((\delta_{1}, \cdots \delta_{k}),\widehat{1}) $ equals to 
$ \lambda$ 
and $\lambda$ is bigger than $t_{1,i}$ for $1 \le i \le N$ in 
the total order on $T^{k} \cup \{ \theta \}$.
Moreover we have 

$c$ is a falling maximal chain 

$\Longleftrightarrow$

$c \in \{ (e, \cdots e) \prec \cdots \prec (e, \cdots e, \delta_{k}) \prec
\cdots (e, \cdots \delta_{k-1},\delta_{k}) \prec \cdots 
\prec (e,e,\delta_{3}, \cdots \delta_{k}) \prec 
\cdots \prec (e,\delta_{2} \cdots \delta_{k})    \ | \ 
{\rm \ each  \ part } \ 
 ( e, \cdots e) \prec \cdots \prec (e, \cdots e, \delta_{k}) ,   
 (e, \cdots e, \delta_{k}) \prec\cdots (e, \cdots \delta_{k-1},\delta_{k})   
 \cdots    
  (e,e,\delta_{3}, \cdots \delta_{k}) \prec \cdots \prec 
  (e,\delta_{2} \cdots \delta_{k}) \  
  {\rm are \ falling \ maximal \ chain}\}$
  
  Hence we have

$\mu(NC_{k}(W) \cup \{\widehat{1}\})$
 
$= \Sigma_
{(e,\delta_{2}, \cdots ,\delta_{k}) \in {\rm maxs}}$

$(-1)^{{\rm rank}(\delta_{2})} \times \{$
the number of falling maximal chains from $e$ to $\delta_{2}$ with respect to 
$\lambda$ $\}$

$\times (-1)^{{\rm rank}(\delta_{3})} \times \{$
the number of falling maximal chains from $e$ to $\delta_{3}$ with respect to 
$\lambda$ $\}$

\begin{center}
$\vdots$
\end{center}

$ \times (-1)^{{\rm rank}(\delta_{k})} \times \{$
the number of falling maximal chains from $e$ to $\delta_{k}$ with respect to 
$\lambda$ $\}$

$= \Sigma_
{(e,\delta_{2}, \cdots ,\delta_{k}) \in {\rm maxs}}$
$\mu([e,\delta_{2}]) \times \mu([e,\delta_{3}]) \cdots \mu([e,\delta_{k}])$

$=\Sigma_
{(\delta_{1}, \cdots ,\delta_{k-1}) : \delta_{1} \cdots \delta_{k-1}=c \ 
{\rm and} l(\delta_{1}) + \cdots l(\delta_{k-1}) = n-1}$
$\mu([e,\delta_{1}]) \times \mu([e,\delta_{2}]) \cdots \mu([e,\delta_{k-1}])$.

Now we have the following proposition.

\begin{prop}

$\Sigma_{(\delta_{1}, \cdots ,\delta_{k-1}) ,  l( \delta_{1} \cdot \delta_{2} 
\cdots \delta_{i}) = l(\delta_{1}) + \cdots l(\delta_{i}), 
\delta_{1} \cdots \delta_{k-1} = c} \mu([e,\delta_{1}]) 
\cdots \mu(e,\delta_{k-1}) = (-1)^{n} {\rm Cat}_{+}^{k-1}(W)$.

\end{prop}

$Proof$

In \cite{armstrong}, Armstrong 
showed that $\mu(NC^{(k)}(W) \cup \{ \widehat{0} \})= 
\mu(NC_{(k)}(W) \cup \{\widehat{1}\})=
(-1)^{n} {\rm Cat}_{+}^{(k-1)}(W)$.
From the view of the 
EL-labeling introduced by Armstrong and Thomas, 
we have 

$c$ is a falling maximal chain 

$\Longleftrightarrow$

$c \in \{ (e, \cdots e) \prec \cdots \prec (e, \cdots e, \delta_{k-1}) \prec
\cdots (e, \cdots \delta_{k-2},\delta_{k-1}) \prec \cdots 
\prec (e,\delta_{2}, \cdots \delta_{k-1}) \prec 
\cdots \prec (\delta_{1} \cdots \delta_{k-1})    \ | \ 
{\rm \ each  \ part } \ 
 ( e, \cdots e) \prec \cdots \prec (e, \cdots e, \delta_{k-1}) ,   
 (e, \cdots e, \delta_{k-1}) \prec\cdots (e, \cdots \delta_{k-2},\delta_{k-1})   \cdots    
  (e,\delta_{2}, \cdots \delta_{k-1}) \prec \cdots \prec 
  (\delta_{1} \cdots \delta_{k-1}) \  
  {\rm are \ falling \ maximal \ chain}\}$.
  
Hence we have 

$(\mu(NC_{(k)}(W) \cup \{\widehat{1}\}))$

$= \Sigma_
{(\delta_{1}, \cdots ,\delta_{k-1}) ; 
( \delta_{1} \cdot \delta_{2} 
\cdots \delta_{i}) = l(\delta_{1}) + \cdots l(\delta_{i}), 
\delta_{1} \cdots \delta_{k-1} = c}$

$(-1)^{{\rm rank}(\delta_{1})} \times \{$
the number of falling maximal chains from $e$ to $\delta_{1}$ with respect to 
$\lambda$ $\}$

$\times (-1)^{{\rm rank}(\delta_{2})} \times \{$
the number of falling maximal chains from $e$ to $\delta_{2}$ with respect to 
$\lambda$ $\}$

\begin{center}
$\vdots$
\end{center}

$ \times (-1)^{{\rm rank}(\delta_{k-1})} \times \{$
the number of falling maximal chains from $e$ to $\delta_{k-1}$ $\}$

$ = \Sigma_{(\delta_{1}, \cdots ,\delta_{k-1}) ,  l( \delta_{1} \cdot \delta_{2} 
\cdots \delta_{i}) = l(\delta_{1}) + \cdots l(\delta_{i}), 
\delta_{1} \cdots \delta_{k-1} = c} \mu([e,\delta_{1}]) 
\cdots \mu(e,\delta_{k-1})$.  

Hence we obtain the derived result.
\ \ \ $\Box$

Now we have 

$\Sigma_{x \in {\rm maxs}} \mu(\widehat{0},x) = 
\Sigma_{(\delta_{1}, \cdots ,\delta_{k}) ,  l( \delta_{1} \cdot \delta_{2} 
\cdots \delta_{i}) = l(\delta_{1}) + \cdots l(\delta_{i}) {\rm for } 1 \le i \le k, 
\delta_{1} \cdots \delta_{k} = c} \mu([e,\delta_{1}]) 
\cdots \mu(e,\delta_{k-1}) = (-1)^{n} {\rm Cat}_{+}^{k-1}(W)
$

This complete the proof of our Theorem {\rmfamily \ref{theorem}}.

\vspace{5mm}

{\large \textbf{Acknowledgement}}

The author wishes to thank Professor Christian Krattenthaler,  
Professor Jun Morita for their 
valuable advices.

\renewcommand{\refname}{REFERENCE}

\end{document}